\documentclass[11pt]{article}%
\usepackage{hyperref}
\usepackage{xcolor}
\usepackage{amsmath}
\usepackage{amsfonts}
\usepackage{amssymb}
\usepackage{graphicx}%
\providecommand{\U}[1]{\protect\rule{.1in}{.1in}}
\newtheorem{theorem}{Theorem}

\newtheorem{proposition}[theorem]{Proposition}

\begin{document}

\title{Asymptotic and exterior Dirichlet problems for the minimal surface equation in
the Heisenberg group with a balanced metric}
\author{Fidelis Bittencourt
\and Edson S. Figueiredo
\and Pedro Fusieger
\and Jaime Ripoll}
\date{}
\maketitle

\begin{abstract}
It is proved that the Heisenberg group $\operatorname*{Nil}\nolimits_{3}$ with
a balanced metric, the sum of the left and right invariant metrics, splits as
a Riemannian product $\mathbb{T\times Z}$, where $\mathbb{T}$ is a totally
geodesic surface and $\mathbb{Z}$ the center of $\operatorname*{Nil}%
\nolimits_{3}.$ It is then proved the existence of complete properly embedded
minimal surfaces in $\operatorname*{Nil}\nolimits_{3}$ by solving the
asymptotic Dirichlet problem for the minimal surface equation on $\mathbb{T}$.
It is also proved the existence of complete properly embedded minimal surfaces
foliating an open set of $\operatorname*{Nil}\nolimits_{3}$ having as boundary
a given curve $\Gamma$ in $\mathbb{T},$ satisfying the exterior circle
condition, by solving the exterior Dirichlet problem for the minimal surface
equation in the unbounded connected component of $\mathbb{T}\backslash\Gamma$.

\end{abstract}

\section{Introduction}

\qquad The study of minimal surfaces in $3-$dimensional Lie groups with a left
invariant metric has recently been attracting the attention of many
mathematicians (\cite{MMP}). In special, in the Heisenberg group%
\[
\operatorname*{Nil}\nolimits_{3}=\left\{  \left(  x,y,z\right)  :=\left(
\begin{array}
[c]{ccc}%
1 & x & z\\
0 & 1 & y\\
0 & 0 & 1
\end{array}
\right)  ,\text{ }x,y,z\in\mathbb{R}\right\}
\]
(see \cite{ADR}, \cite{DL}, \cite{YZ}, \cite{FMP}, \cite{Rit}, \cite{IJL}).

Since a left invariant metric on a Lie group is completely determined by its
value on the Lie algebra of the group and the group operation, the Riemannian
geometry of the group is closely related to the Lie group structure. To bring
up this connection, between the group and the Riemannian structures, is one of
the motivations for studying minimal surfaces in a Lie group with a left
invariant metric. In this paper we consider a \emph{balanced} \emph{metric} on
$\operatorname*{Nil}\nolimits_{3}$, namely:%
\begin{equation}
\left\langle u,v\right\rangle _{g}=\left\langle d\left(  L_{g}\right)
_{g}^{-1}u,d\left(  L_{g}\right)  _{g}^{-1}v\right\rangle _{e}+\left\langle
d\left(  R_{g}\right)  _{g}^{-1}u,d\left(  R_{g}\right)  _{g}^{-1}%
v\right\rangle _{e}, \label{bm}%
\end{equation}
where $e$ is the neutral element of $\operatorname*{Nil}\nolimits_{3}$ (the
identity matrix), the inner product of $T_{e}\operatorname*{Nil}\nolimits_{3}$
given by%
\[
\left\langle \left(
\begin{array}
[c]{ccc}%
0 & a_{1} & c_{1}\\
0 & 0 & b_{1}\\
0 & 0 & 0
\end{array}
\right)  ,\left(
\begin{array}
[c]{ccc}%
0 & a_{2} & c_{2}\\
0 & 0 & b_{2}\\
0 & 0 & 0
\end{array}
\right)  \right\rangle _{e}=a_{1}a_{2}+b_{1}b_{2}+c_{1}c_{2},
\]
$g\in\operatorname*{Nil}\nolimits_{3},$ $u,v\in T_{g}\operatorname*{Nil}%
\nolimits_{3}$ and $L_{g},R_{g}:\operatorname*{Nil}\nolimits_{3}%
\rightarrow\operatorname*{Nil}\nolimits_{3}$ are the left and right
translations, $L_{g}(x)=gx,$ $R_{g}(x)=xg.$

Clearly, the metric (\ref{bm}) is completely determined by the its value on
the Lie algebra of the group and by the group operation and, thus, the
geometry of $\operatorname*{Nil}\nolimits_{3}$ with the balanced metric is
expected to be closely related to the group structure of $\operatorname*{Nil}%
\nolimits_{3}$. Indeed, we prove that with the balanced metric
$\operatorname*{Nil}\nolimits_{3}$ splits as a Riemannian product of its
center and a totally geodesic surface$.$ Actually, the geometry of
$\operatorname*{Nil}\nolimits_{3}$ with the balanced metric presents
surprising properties, as shown in our first result:

\begin{theorem}
\label{gnil}Setting
\[
\mathbb{T}=\left\{  \left(  x,y,\frac{xy}{2}\right)  \in\operatorname*{Nil}%
\nolimits_{3},\text{ }x,y\in\mathbb{R}\right\}
\]
we have:

\begin{itemize}
\item[(a)] $\mathbb{T}$ is a totally geodesic surface of $\operatorname*{Nil}%
\nolimits_{3}.$

\item[(b)] $\left(  \operatorname*{Nil}\nolimits_{3},\left\langle \text{ ,
}\right\rangle \right)  $ is isometric to the Riemannian product
$\mathbb{T}\times\mathbb{Z},$ where $\mathbb{Z}=\left\{  \left(  0,0,z\right)
\text{
$\vert$
}z\in\mathbb{R}\right\}  \ $is the center of $\operatorname*{Nil}%
\nolimits_{3}$

\item[(c)] $\mathbb{S}^{1}=\left\{  e^{i\theta}\right\}  _{\theta\in
\mathbb{R}}$ acts isometrically on $\operatorname*{Nil}\nolimits_{3}$ by%
\begin{gather*}
e^{i\theta}\left(  x,y,z\right)  =\\
\left(  x\cos\theta-y\sin\theta,x\sin\theta+y\cos\theta,\frac{\left(
x\sin\theta+y\cos\theta\right)  \left(  x\cos\theta-y\sin\theta\right)  }%
{2}+z\right)  .
\end{gather*}
This action leaves $\mathbb{T}$ invariant and pointwise fixed the center of
$\operatorname*{Nil}\nolimits_{3}.$

\item[(d)] The arc length geodesics of $\mathbb{T}$ passing through the
identity $e$ are%
\[
\gamma_{\theta}(t)=\left(
\begin{array}
[c]{ccc}%
1 & \frac{t}{2}\left(  \cos\theta-\sin\theta\right)   & \frac{\left(
\cos\theta-\sin\theta\right)  \left(  \sin\theta+\cos\theta\right)  t^{2}}%
{8}\\
0 & 1 & \frac{t}{2}\left(  \sin\theta+\cos\theta\right)  \\
0 & 0 & 1
\end{array}
\right)  ,\text{ }t\in\mathbb{R},\text{ }\theta\in\left[  0,2\pi\right)  .
\]
Being $r(p)=d(p,e),$ $d$ the Riemannian distance in $\mathbb{T},$ we have%
\begin{equation}
r\left(  x,y,z\right)  =\sqrt{2}\sqrt{x^{2}+y^{2}}.\label{erre}%
\end{equation}

\item[(e)] The sectional curvature $K$ of $\mathbb{T}$ is
\[
K_{p}=-4\frac{r^{2}+12}{\left(  r^{2}+8\right)  ^{2}},\text{ \ }p\in
\mathbb{T},
\]
where $r=r(p).$
\end{itemize}
\end{theorem}

Notice that from Theorem \ref{gnil}, item (e), it is clear that with the
balanced metric $\operatorname*{Nil}\nolimits_{3}$ is not a homogeneous manifold.

The Riemannian splitting $\operatorname*{Nil}\nolimits_{3}\simeq
\mathbb{T}\times\mathbb{Z}$ allows the construction of complete properly
embedded minimal surfaces of $\operatorname*{Nil}\nolimits_{3}$ with the
balanced metric by solving the asymptotic Dirichlet problem in $\mathbb{T}$
(see \cite{RT} for details on the this problem). Indeed, it is an immediate
consequence of Corollary 1.2 of \cite{RTe}:

\begin{theorem}
\label{ap}The asymptotic Dirichlet problem
\[
\left\{
\begin{array}
[c]{l}%
\mathcal{M}\left[  u\right]  :=\operatorname{div}\frac{\nabla u}%
{\sqrt{1+\left\Vert \nabla u\right\Vert ^{2}}}=0\text{ in }\mathbb{T}\\
u|_{\partial_{\infty}\mathbb{T}}=\varphi
\end{array}
\right.
\]
has one and only one solution $u\in C^{\infty}\left(  \mathbb{T}\right)  \cap
C^{0}\left(  \overline{\mathbb{T}}\right)  $ for any $\varphi\in C^{0}\left(
\partial_{\infty}\mathbb{T}\right)  .$
\end{theorem}

Notice that the graph $G(u)=\left\{  \left(  p,u(p)\right)  \text{
$\vert$
}p\in\mathbb{T}\right\}  $ is a complete minimal surface in
$\operatorname*{Nil}_{3}$ with respect to the balanced metric (\ref{bm}). In
matricial terms, using the identification $u(x,y)=u(x,y,(xy)/2),$ the graph of
$u$ is%
\[
G(u)=\left\{  \left(
\begin{array}
[c]{ccc}%
1 & x & \frac{xy}{2}+u(x,y)\\
0 & 1 & y\\
0 & 0 & 1
\end{array}
\right)  ,\text{ \ \ }x,y\in\mathbb{R}\right\}  .
\]
In particular, it follows from Theorem \ref{ap} the existence of an infinite
number of non congruent foliations of $\operatorname*{Nil}\nolimits_{3}$ by
complete properly embedded minimal surfaces transversal to the center of
$\operatorname*{Nil}\nolimits_{3}$. We should remark that Theorem \ref{ap} is
a special case of a family of asymptotic Dirichlet problems which is being
investigated on a big class of PDE's and in general Hadamard manifolds. This
is an active area of research which produced a vast numbers of papers in
recent years (see \cite{RT}).

We also study the exterior Dirichlet problem (EDP) for the minimal surface
equation in $\mathbb{T}$ namely, the existence and uniqueness of solutions to
the PDE $\mathcal{M}\left[  u\right]  =0$ on a domain $\Lambda\subset
\mathbb{T}$, with a prescribed data at $\partial\Lambda,$ when $\mathbb{T}%
\backslash\overline{\Lambda}$ is a bounded domain. This problem can of course
be considered in any complete non compact Riemannian manifold, and has an old
history about which it would be interesting to say some words. It was first
studied in the Euclidean space by J. C. C. Nitsche who proved that possible
solutions have at most linear growth (\cite{N}). R. Osserman proved the
existence of a boundary data at $\partial\Lambda$, even when $\mathbb{R}%
^{2}\backslash\Lambda$ is a disk, for which the EDP has no solution
(\cite{O}). R. Krust proved that the solutions of the EDP having the same
Gauss map at infinity form a foliation if there are at least two solutions
(\cite{Kr}). Krust's result was improved by J. Ripoll and F. Tomi in
\cite{RT1} where they proved the existence of a minimal and of a maximal
solutions and also the existence of a boundary data admitting exactly one
solution. The results of Krust were extended to arbitrary dimensions by E.
Kuwert in \cite{K}. Still in the Euclidean space the EDP was also studied in
\cite{KT}, \cite{R} and, in the Riemannian setting, in \cite{ER}. Our result
on the EDP in $\operatorname*{Nil}\nolimits_{3}$ is based on \cite{R} and
\cite{ER}:

\begin{theorem}
\label{ep}Let $\Omega$ be a $C^{2,\alpha}$ domain of $\mathbb{T}$\ satisfying
the exterior geodesic circle condition, namely: given $p\in\Omega$, there
exists a geodesic circle passing through $p$ which is the boundary of a
geodesic disk containing $\Omega.$ Set $\Lambda=\mathbb{T}\backslash\Omega$.
Then, for each $s\geq0$ the exterior Dirichlet problem in $\Lambda$:%
\[
\left\{
\begin{array}
[c]{l}%
\mathcal{M}\left[  u\right]  =0\text{ in }\Lambda\\
u|_{\partial\Lambda}=0
\end{array}
\right.
\]
has a solution $u_{s}\in C^{2,\alpha}\left(  \overline{\Lambda}\right)  $ such
that $\sup_{\partial\Lambda}\left\Vert \nabla u\right\Vert =s.$ Moreover, the
graphs of $u_{s}$ form a foliation of an open subset of $\operatorname*{Nil}%
\nolimits_{3}$ and
\[
\limsup_{r(p)\rightarrow\infty}\left\vert u_{s_{1}}(p)-u_{s_{2}}(p)\right\vert
>0
\]
if $s_{1}\neq s_{2}.$
\end{theorem}

\section{Preliminaries\label{P}}

\qquad(i) We shall use the following parametrization $\varphi:\mathbb{R}%
^{3}\rightarrow\operatorname*{Nil}\nolimits_{3}$ of $\operatorname*{Nil}%
\nolimits_{3},$
\[
\varphi(x,y,z)=\left(  x,y,\frac{xy}{2}+z\right)
\]
and the corresponding coordinate vector fields%

\[
X=\frac{\partial\varphi}{\partial x},\text{ }Y=\frac{\partial\varphi}{\partial
y},\text{ }Z=\frac{\partial\varphi}{\partial z}.
\]

The coefficients\ at $\varphi(x,y,z)$ on the basis $\{X,Y,Z\}$ of the metric
(\ref{bm}) are
\begin{align*}
\left\langle X,X\right\rangle  &  =2+\frac{1}{2}y^{2},~\left\langle
Y,Y\right\rangle =2+\frac{1}{2}x^{2},\\
\left\langle X,Y\right\rangle  &  =-\frac{1}{2}xy,~\left\langle
X,Z\right\rangle =0,\\
\left\langle Y,Z\right\rangle  &  =0,~\left\langle Z,Z\right\rangle =2
\end{align*}
and the Riemannian connection $\nabla$ is given, at $\ \varphi(x,y,z),$ by%
\begin{align*}
\nabla_{X}X  &  =-\frac{xy^{2}}{2x^{2}+2y^{2}+8}X-\frac{4y+y^{3}}%
{2x^{2}+2y^{2}+8}Y,\\
\nabla_{Y}X  &  =\frac{2y+x^{2}y}{2x^{2}+2y^{2}+8}X+\frac{2x+xy^{2}}%
{2x^{2}+2y^{2}+8}Y,\\
\nabla_{X}Y  &  =\frac{2y+x^{2}y}{2x^{2}+2y^{2}+8}X+\frac{2x+xy^{2}}%
{2x^{2}+2y^{2}+8}Y,\\
\nabla_{Y}Y  &  =-\frac{4x+x^{3}}{2x^{2}+2y^{2}+8}X-\frac{x^{2}y}%
{2x^{2}+2y^{2}+8}Y
\end{align*}
and%
\[
\nabla_{Y}Z=\nabla_{Z}Y=\nabla_{Z}X=\nabla_{X}Z=\nabla_{Z}Z=0.
\]

(ii) One may see that the curve $\gamma:\mathbb{R\rightarrow T},$
$\gamma(t)=\varphi\left(  t/2,t/2,0\right)  $ is an arc length geodesic of
$\operatorname*{Nil}\nolimits_{3}$. Indeed, note that%

\[
\gamma^{\prime}(t)=\frac{1}{2}\left(  X+Y\right)  .
\]
and hence%

\begin{align*}
4\nabla_{\gamma^{\prime}}\gamma^{\prime}(t)  &  =\nabla_{X+Y}\left(
X+Y\right) \\
&  =\left(  -\frac{xy^{2}}{2x^{2}+2y^{2}+8}X-\frac{4y+y^{3}}{2x^{2}+2y^{2}%
+8}Y\right) \\
&  +2\left(  \frac{2y+x^{2}y}{2x^{2}+2y^{2}+8}X+\frac{2x+xy^{2}}{2x^{2}%
+2y^{2}+8}Y\right) \\
&  +\left(  -\frac{4x+x^{3}}{2x^{2}+2y^{2}+8}X-\frac{x^{2}y}{2x^{2}+2y^{2}%
+8}Y\right)  =0
\end{align*}
since $x=x(t)=t/2=y(t)=y.$ Moreover%
\begin{align*}
\left\Vert \gamma^{\prime}\right\Vert ^{2}  &  =\left(  x^{\prime}\right)
^{2}\left(  2+\frac{y^{2}}{2}\right)  +\left(  y^{\prime}\right)  ^{2}\left(
2+\frac{x^{2}}{2}\right)  +2x^{\prime}y^{\prime}\left(  -\frac{xy}{2}\right)
\\
&  =\frac{1}{4}\left(  \left(  2+\frac{y^{2}}{2}\right)  +\left(
2+\frac{x^{2}}{2}\right)  -xy\right)  =1.
\end{align*}

\section{Proof of Theorem \ref{gnil}}

\qquad(a) $\mathbb{T}$ is a totally geodesic surface of $\operatorname*{Nil}%
\nolimits_{3}$ since one may see that $Z$ is a Killing field orthogonal to
$\mathbb{T}$.

(b) An isometry $\Psi:\mathbb{T}\times\mathbb{Z\rightarrow}\operatorname*{Nil}%
\nolimits_{3}$ is given by%
\[
\Psi\left(  \left(  \left(
\begin{array}
[c]{ccc}%
1 & x & \frac{1}{2}xy\\
0 & 1 & y\\
0 & 0 & 1
\end{array}
\right)  ,\left(
\begin{array}
[c]{ccc}%
1 & 0 & t\\
0 & 1 & 0\\
0 & 0 & 1
\end{array}
\right)  \right)  \right)  =\left(
\begin{array}
[c]{ccc}%
1 & x & \frac{1}{2}xy+t\\
0 & 1 & y\\
0 & 0 & 1
\end{array}
\right)  .
\]

(c) Follows from direct computations

(d) The first part follows from Section \ref{P} (i) and (c). For the remaining
part, given $p\in\mathbb{T}$, $p=\varphi(x,y,0)$, taking $t=\sqrt{2\left(
x^{2}+y^{2}\right)  }$ there is $\theta\in\mathbb{R}$ such that
\[
x=\frac{t}{2}\left(  \cos\theta-\sin\theta\right)  ,~y=\frac{t}{2}\left(
\sin\theta+\cos\theta\right)
\]
and hence%
\[
x^{2}+y^{2}=\left(  \frac{t}{2}\left(  \cos\theta-\sin\theta\right)  \right)
^{2}+\left(  \frac{t}{2}\left(  \sin\theta+\cos\theta\right)  \right)
^{2}=\frac{1}{2}t^{2}.
\]
Since $t$ is the arc length of the geodesic
\[
t\rightarrow\left(  \frac{t}{2}\left(  \cos\theta-\sin\theta\right)
,y=\frac{t}{2}\left(  \sin\theta+\cos\theta\right)  ,\frac{\left(  \cos
\theta-\sin\theta\right)  \left(  \sin\theta+\cos\theta\right)  t^{2}}%
{8}\right)
\]
we obtain (\ref{erre}).

(e) Note that the vector fields $X$ and $Y$ are tangent to $\mathbb{T}$ and
$\left[  X,Y\right]  =0.$ Hence%
\[
\sqrt{\left\Vert X\right\Vert ^{2}\left\Vert Y\right\Vert ^{2}-\left\langle
X,Y\right\rangle ^{2}}K=\left\langle R(X,Y)X,Y\right\rangle =\left\langle
\nabla_{Y}\nabla_{X}X-\nabla_{X}\nabla_{Y}X,Y\right\rangle .
\]
Using the computations done in the preliminary section, after somewhat long
but straightforward calculations one arrives to
\[
K=-\frac{2\left(  x^{2}+y^{2}+6\right)  }{\left(  x^{2}+y^{2}+4\right)  ^{2}%
}=-4\frac{r^{2}+12}{\left(  r^{2}+8\right)  ^{2}}%
\]
$\allowbreak$concluding the proof of Theorem \ref{gnil}.

\section{The exterior Dirichlet problem}

\qquad The catenoids, besides being interesting on their own, pieces of them
provide explicit examples of solutions of the exterior Dirichlet problem
foliating open subsets of $\mathbb{T\times R}$. From straightforward
calculations, one has:

\begin{proposition}
Let $c>0,$ $t_{0}\geq\sqrt{\sqrt{c^{2}+16}-4}$ be given. Then the rotation of
the curve%
\begin{equation}
\alpha(t)=\left(  \frac{t}{2},\frac{t}{2},\frac{c}{\sqrt{2}}\int_{t_{0}}%
^{t}\frac{ds}{\sqrt{s^{2}\left(  s^{2}+8\right)  -c^{2}}}\right)  ,\text{
\ }t\in\mathbb{R}, \label{cat}%
\end{equation}
around around the $z-$ axis is a catenoid that is, a complete rotationally
invariant minimal surface in $\operatorname*{Nil}\nolimits_{3}$.
\end{proposition}

\subsection{Proof of Theorem \ref{ep}}

\qquad We begin by introducing some notation and remarking some facts to be
used in the course of the proof.

Let $r:\Lambda\rightarrow\left[  0,\infty\right)  $ be the distance function
to $\partial\Omega,$ $r(x)=d(x,\partial\Omega),$ where $d$ is the Riemannian
distance in $\mathbb{T}$. Given $r,$ set%
\[
\Omega_{r}=\left\{  x\in\Lambda\text{
$\vert$
}r(x)<r\right\}  .
\]
One may see that $\Omega_{r}$ satisfies the exterior geodesic circle
condition. Then, given $p\in\partial\Omega_{r}\backslash\partial\Omega,$ there
is a geodesic circle enclosing a geodesic disk $D^{p}\subset\mathbb{T}$
passing through $p$ and containing $\Omega_{r}.$ By the tangency principle,
the curvature of this geodesic circle is smaller than or equal to the
curvature of $\partial\Omega_{r}\backslash\partial\Omega$ at $p.$ Take
$R(r)>0$ such that%
\[
\cup_{p\in\partial\Omega_{r}\backslash\partial\Omega}D^{p}\subset D_{R(r)},
\]
where $D_{R(r)}$ is the geodesic disk of $\mathbb{T}$ centered at $e$ with
radius $R(r).$ By Theorem \ref{gnil} (e), the maximum of the sectional
curvature $K_{r}$ of $D_{R(r)}$ is%
\[
K_{r}=-\frac{4\left(  R(r)^{2}+12\right)  }{\left(  R(r)^{2}+8\right)  ^{2}}.
\]
By the Hessian comparison theorem the curvature of the geodesic circles of
$\mathbb{T}$ contained in $D_{R(r)},$ oriented inwards, are greater than or
equal to the curvature of the geodesic circles of the hyperbolic plane
$\mathbb{H}^{2}\left(  K_{r}\right)  $ of sectional curvature $K_{r}.$ In
turn, the curvature $k_{r}$ of any geodesic circles of $\mathbb{H}^{2}\left(
K_{r}\right)  $ has the lower bound%
\[
k_{r}>\sqrt{\frac{4\left(  R(r)^{2}+12\right)  }{\left(  R(r)^{2}+8\right)
^{2}}}>\frac{2R(r)}{R(r)^{2}+8}.
\]
Since $\partial\Omega$ is compact it follows from the triangle inequality
that
\[
\lim_{r\rightarrow\infty}\frac{R(r)}{r}=1
\]
and thus, there is $\alpha>0$ such that $r-\alpha\leq R(r)\leq r+\alpha$ for
all $r\geq0.$ It follows that%
\begin{equation}
k_{r}>\frac{2\left(  r-\alpha\right)  }{\left(  r+\alpha\right)  ^{2}+8}.
\label{bb}%
\end{equation}

For the given $s\geq0$ we now find a subsolution $v_{s}\in C^{2}\left(
\overline{\Lambda}\right)  $ of
\[
\mathcal{M}\left[  u\right]  :=\operatorname{div}\left(  \frac{\nabla u}%
{\sqrt{1+\left\Vert \nabla u\right\Vert ^{2}}}\right)
\]
of the form $v_{s}(x)=f(r(x)),$ $x\in\Lambda,$ for some $f\in C^{2}\left(
\left[  0,\infty\right)  \right)  ,$ such that
\begin{align}
v_{s}|_{\partial\Omega}  &  =0\nonumber\\
\left\Vert \nabla v_{s}(x)\right\Vert  &  =s,\text{ for any }x\in
\partial\Omega\label{cond}\\
\sup_{\Lambda}\left\Vert v_{s}\right\Vert  &  <\infty.\nonumber
\end{align}

We have%
\[
\mathcal{M}\left[  u\right]  =\operatorname{div}\left(  \frac{\nabla v_{s}%
}{\sqrt{1+\left\Vert \nabla v_{s}\right\Vert ^{2}}}\right)  =\frac
{f^{\prime\prime}(r)+f^{\prime}(r)\left(  1+f^{\prime2}(r)\right)  \Delta
v_{s}}{\left(  1+f^{\prime2}(r)\right)  ^{\frac{3}{2}}}
\]
where $\Delta$ is the Laplacian operator. Since $\Delta v_{s}$ is the
curvature the $\partial\Omega_{r}\backslash\partial\Omega,$ assuming that
$f^{\prime}(r)\geq0,$ it follows from (\ref{bb}) that $\mathcal{M}\left[
u\right]  \geq0$ if%

\[
f^{\prime\prime}(r)+f^{\prime}(r)\left(  1+f^{\prime2}(r)\right)  f^{\prime
}(r)\frac{2\left(  r-\alpha\right)  }{\left(  r+\alpha\right)  ^{2}+8}\geq
f^{\prime\prime}(r)+f^{\prime}(r)\frac{2\left(  r-\alpha\right)  }{\left(
r+\alpha\right)  ^{2}+8}\geq0.
\]

Solving the last inequality to have an equality, we obtain%

\[
f(r)=\int_{0}^{r}\frac{ce^{\sqrt{2}\alpha\arctan\left(  \frac{t+\alpha
}{2^{3/2}}\right)  }}{\left(  t+\alpha\right)  ^{2}+8}dt
\]
where $c$ is any constant. Clearly $f(0)=0$ and%

\[
f^{\prime}(0)=\frac{c\exp\left(  \sqrt{2}\alpha\arctan\frac{1}{4}\sqrt
{2}\alpha\right)  }{\alpha^{2}+8}%
\]
so that we can choose $c$ such that $f^{\prime}(0)=s.$ Since
\[
0\leq e^{\sqrt{2}\alpha\arctan\left(  \frac{t+\alpha}{2^{3/2}}\right)  }\leq
e^{\sqrt{2}\alpha\frac{\pi}{2}},\text{ }t\geq0,
\]
it follows that $f$ is bounded. Hence, $v_{s}(x)=f(r(x)),$ $x\in
\overline{\Lambda},$ is a subsolution of $\mathcal{M}$ satisfying conditions
(\ref{cond}).

Given $m\in\mathbb{N}$, set
\begin{align*}
T_{m}  &  =\left\{  t\geq0\text{
$\vert$
}\exists u_{t}\in C^{2,\alpha}\left(  \overline{\Omega}_{m}\right)  \text{
s.t. }\mathcal{M}\left[  u_{t}\right]  =0,\text{ }\right. \\
&  \left.  u_{t}|_{\partial\Omega}=0,\text{ }u_{t}|_{\partial\Omega
_{m}\backslash\partial\Omega}=t,\text{ }\sup\nolimits_{\partial\Omega
}\left\Vert \nabla u_{t}\right\Vert \leq s\right\}  .
\end{align*}

We have $T_{m}\neq\emptyset,$ since $0\in T_{m}$. Give $\varepsilon>0,$ we
prove that $t<v_{s+\varepsilon}|_{\partial\Omega_{m}\backslash\partial\Omega
}.$ By contradiction, assume that $t\geq v_{s+\varepsilon}|_{\partial
\Omega_{m}\backslash\partial\Omega}.$ Since
\[
\left\Vert \nabla v_{s+\varepsilon}\right\Vert _{\partial\Omega}%
=\inf_{\partial\Omega}\left\Vert \nabla v_{s+\varepsilon}\right\Vert
_{\partial\Omega}=s+\varepsilon>\sup_{\partial\Omega}\left\Vert \nabla
u_{t}\right\Vert _{\partial\Omega}%
\]
there is a neighborhood $U$ of $\partial\Omega$ in $\Omega_{m}$ such that
$u_{t}(x)<v_{s+\varepsilon}(x)$ for all $x\in U\backslash\partial\Omega_{m}.$
Since $t<v_{s+\varepsilon}|_{\partial\Omega_{m}\backslash\partial\Omega}$
there exists a domain $V\subset\Omega_{m}$ containing $U$ such that
$u_{t}|_{\partial V}=v_{s+\varepsilon}|_{\partial V},$ what is an absurd since
$v_{s+\varepsilon}|_{V}$ is a subsolution of $\mathcal{M}\ $in $V$ coinciding
with $u_{t}$ at $\partial V$ (this follows from the comparison principle,
Proposition 3.1 of \cite{RT})$.$ Letting $\varepsilon\rightarrow0$ we have
$t\leq v_{s}|_{\partial\Omega_{m}\backslash\partial\Omega}.$ It follows that
$T_{m}$ is bounded and we may set%
\[
t_{m}=\sup T_{m}<\infty.
\]
We prove that $t_{m}\in T_{m},$ by first proving that there is a constant $C$,
not depending on $m,$ such that if $t\in T_{m},$ then $\sup_{\partial
\Omega_{m}}\left\Vert \nabla u_{t}\right\Vert _{\partial\Omega}\leq C.$ Let
$t\in T_{m}$ be given$.$ By definition of $T_{m}$ we have $\sup_{\partial
\Omega}\left\Vert \nabla u_{t}\right\Vert _{\partial\Omega}\leq s.$

Setting $z_{m}=v_{s}|_{\partial\Omega_{m}\backslash\partial\Omega}$ we have,
as proved above, $z_{m}\geq t$. Moreover, the function $w_{s}:=v_{s}-\left(
z_{m}-t\right)  $ is a subsolution of $\mathcal{M}$ in $\Omega_{m}$ and%
\begin{align*}
w_{s}|_{\partial\Omega}  &  =-\left(  z_{m}-t\right)  \leq0=u_{t}|_{\Omega}\\
w_{s}|_{\partial\Omega_{m}\backslash\partial\Omega}  &  =t=u_{t}%
|_{\partial\Omega_{m}\backslash\partial\Omega}.
\end{align*}

By the comparison principle we obtain%
\begin{align*}
w_{s}  &  \leq u_{t}\leq t\\
w_{s}|_{\partial\Omega_{m}\backslash\partial\Omega}  &  =t=u_{t}%
|_{\partial\Omega_{m}\backslash\partial\Omega}%
\end{align*}
from what we conclude that%
\begin{align*}
\sup_{\partial\Omega_{m}\backslash\partial\Omega}\left\Vert \nabla
u_{t}\right\Vert  &  \leq\sup_{\partial\Omega_{m}\backslash\partial\Omega
}\left\Vert \nabla w_{s}\right\Vert _{\partial\Omega}=f^{\prime}(m)\\
&  =\frac{ce^{\sqrt{2}\alpha\arctan\left(  \frac{m+\alpha}{2^{3/2}}\right)  }%
}{\left(  m+\alpha\right)  ^{2}+8}%
\end{align*}
where
\[
c=\frac{s\left(  \alpha^{2}+8\right)  }{\exp\left(  \sqrt{2}\alpha\arctan
\frac{1}{4}\sqrt{2}\alpha\right)  }.
\]
Since
\[
\sup_{\Omega_{m}}\left\vert u_{t}\right\vert \leq t_{m}\leq\sup_{\Lambda
}\left\vert v_{s}\right\vert <\infty
\]
for any $t\in T_{m}$ we may then conclude the $C^{1}$ norm of any solution
$u_{t}$ for $t\in T_{m}$ has an a-priori bound that depends only on $s$ and
$\alpha.$ Hence, from elliptic PDE linear theory there is $C=C(s,\alpha)$ such
that $\left\Vert u_{t}\right\Vert _{C^{2,\alpha}\left(  \overline{\Omega
}\right)  }\leq C$ for any $t\in T_{m}$ (see Ch. 2.1 of \cite{RT}).

Consider now a sequence $\left\{  z_{j}\right\}  \subset T_{m}$ converging to
$t_{m}$ as $j$ goes to infinity. For each $j,$ there is a function $u_{j}\in
C^{2,\alpha}(\overline{\Omega}_{m})$ such that $\mathcal{M}\left(
u_{j}\right)  =0,$ $u_{j}|_{\partial\Omega}=0$ and $u_{j}|_{\Gamma_{m}}%
=z_{j}.$ Since $\left\Vert u_{j}\right\Vert _{C^{2,\alpha}\left(
\overline{\Omega}\right)  }\leq C$ there is a subsequence of $\left\{
u_{j}\right\}  $ that converges uniformly in $\overline{\Omega}_{m}$ on the
$C^{2}$ norm to a solution $w_{m}\in C^{2}(\overline{\Omega_{m}})$ of
$\mathcal{M}=0$ in $\Omega_{m}$. From PDE regularity $w_{m}\in C^{2,\alpha
}(\overline{\Omega}_{m})$ (\cite{GT})$.$

The function $w_{m}$ is a solution of $\mathcal{M}=0$ in $\Omega_{m}$ that
satisfies $w_{m}|_{\partial\Omega}=0$, $w_{m}|_{\Gamma_{m}}=t_{m}$ and
$\sup_{\partial\Omega_{m}}\left\vert \operatorname{grad}w_{m}\right\vert \leq
s.$ It follows that $t_{m}\in T_{m},$ that is, $w_{m}=u_{t_{m}}$. Moreover, it
follows from the implicit function theorem $\sup_{\partial\Omega_{m}%
}\left\vert \operatorname{grad}w_{m}\right\vert <s$ leads to a contradiction
so that $\sup_{\partial\Omega}\left\vert \operatorname{grad}w_{m}\right\vert
=s.$

Since the bound $\left\Vert w_{m}\right\Vert _{C^{2,\alpha}\left(
\overline{\Lambda}\right)  }\leq C$ does not depend on $m$ there is a
subsequence of $\{w_{m}\}$ converging uniformly $C^{2}$ on compact subsets of
$\overline{\Lambda}$ to a solution $u_{s}\in C^{2,\alpha}(\overline{\Lambda})$
of $\mathcal{M}=0$ in $\Omega$ satisfying $u_{s}|_{\partial\Omega}=0$ and
$\sup_{\partial\Omega}\left\vert \operatorname{grad}u_{s}\right\vert =s.$

We prove that given if $0\leq s_{1}<s_{2}\ $then $u_{s_{1}}<u_{s_{2}}$ on
$\Lambda.$ Suppose that $u_{s_{i}}=\lim_{m}w_{m}^{i},$ where $w_{m}^{i}\in
C^{2,\alpha}\left(  \overline{\Omega}_{m}\right)  $ is a solution as above
$i=1,2$. Given $p\in\Lambda,$ there is $m\in\mathbb{N}$ such that $p\in
\Omega_{m}.$ It is clear that $w_{n}^{1}|_{\partial\Omega_{n}\backslash
\partial\Omega}<w_{n}^{2}|_{\partial\Omega_{n}\backslash\partial\Omega}$ and
then, from the comparison principle, $w_{n}^{1}(p)<w_{n}^{2}(p)$ for all
$n\geq m.$ It follows that $u_{s_{1}}(p)\leq u_{s_{2}}(p)$ for all
$p\in\Lambda.$ From the maximum principle, $u_{s_{1}}(p)<u_{s_{2}}(p)$ for all
$p\in\Lambda.$ This proves that the graphs of $u_{s}$ form a foliation of an
open set of $\mathbb{T\times R}$. We also have
\[
\limsup_{r(p)\rightarrow\infty}\left\vert u_{s_{2}}(p)-u_{s_{1}}(p)\right\vert
=\limsup_{r(p)\rightarrow\infty}\left(  u_{s_{2}}(p)-u_{s_{1}}(p)\right)
\geq0.
\]
Assume that $\limsup_{r(p)\rightarrow\infty}\left(  u_{s_{2}}(p)-u_{s_{1}%
}(p)\right)  =0.$ We claim that for any $\varepsilon>0$, $u_{s_{2}%
}-\varepsilon<u_{s_{1}}.$ Indeed, if
\[
\left\{  p\in\Lambda\text{
$\vert$
}u_{s_{2}}(p)-\varepsilon\geq u_{s_{1}}(p)\right\}
\]
is nonempty for some $\varepsilon>0$ then%
\[
U:=\left\{  p\in\Lambda\text{
$\vert$
}u_{s_{2}}(p)-\frac{\varepsilon}{2}>u_{s_{1}}(p)\right\}
\]
is also nonempty and open. Since $u_{s_{2}}-\varepsilon/2$ is also a solution
of $\mathcal{M}\left[  u\right]  =0,$ the comparison principle implies that
$U$ is not bounded$.$ Then there is a divergence sequence $p_{n}$ in
$\mathbb{T}$ such that $u_{s_{2}}(p_{n})-\varepsilon/2>u_{s_{1}}(p_{n})$ so
that $\limsup_{n\rightarrow\infty}\left(  u_{s_{2}}(p_{n})-u_{s_{1}%
}(p)\right)  \geq\varepsilon/2,$ contradiction. Letting $\varepsilon$ go to
zero in $u_{s_{2}}-\varepsilon<u_{s_{1}}$ we obtain a contradiction with
$u_{s_{1}}<u_{s_{2}}$ on $\Lambda.$ This concludes the proof of the theorem.

e-mail address: \href{mailto:fidelisb@gmail.com}{fidelisb@gmail.com}\\
{\tiny{DEPARTAMENTO DE MATEM\'{A}TICA, UNIVERSIDADE FEDERAL DE SANTA MARIA, SANTA MARIA, BRAZIL}} \\ 

e-mail address: \href{mailto:esidney@gmail.com}{esidney@gmail.com}\\
{\tiny{DEPARTAMENTO DE MATEM\'{A}TICA, UNIVERSIDADE FEDERAL DE SANTA MARIA, SANTA MARIA, BRAZIL}} \\

e-mail address: \href{mailto:fuspedro@gmail.com}{fuspedro@gmail.com}\\
{\tiny{DEPARTAMENTO DE MATEM\'{A}TICA, UNIVERSIDADE FEDERAL DE SANTA MARIA, SANTA MARIA, BRAZIL}} \\

e-mail address: \href{mailto:jaime.ripoll@ufrgs.com}{jaime.ripoll@ufrgs.com}\\
{\tiny{INSTITUTO DE MATEM\'{A}TICA E ESTAT\'{I}STICA, UNIVERSIDADE FEDERAL DO RIO GRANDE DO SUL, PORTO ALEGRE, BRAZIL}}

\end{document}